\documentclass[24pt]{amsart}
\usepackage{amsmath, amsthm, amscd, amsfonts}

\newtheorem{thm}{Theorem}[section]

\newtheorem{prop}[thm]{Proposition}

\theoremstyle{definition}

\theoremstyle{remark}

\numberwithin{equation}{section}
\newfont{\kh}{msbm10}

\begin{document}

\title{On a Result of Hardy and Ramanujan}

\author{M. Avalin Charsooghi, Y. Azizi, M. Hassani and L. Mollazadeh-Beidokhti}

\address{Mohammad Avalin Charsooghi, Yousof Azizi and Laleh Mollazadeh-Beidokhti, \newline Department of Physics,
Institute for Advanced Studies in Basic Sciences, P.O. Box
45195-1159, Zanjan, Iran}

\email{<avalinch, azizi, laleh>@iasbs.ac.ir}

\address{Mehdi Hassani, \newline Department of Mathematics, Institute for Advanced Studies in Basic
Sciences, P.O. Box 45195-1159, Zanjan, Iran}

\email{mmhassany@member.ams.org}



\begin{abstract}
In this paper, we introduce some explicit approximations for the
summation $\sum_{k\leq n}\Omega(k)$, where $\Omega(k)$ is the total
number of prime factors of $k$.
\end{abstract}

\maketitle

\section{Introduction}
Let $\Omega(k)$ be the total number of prime factors of $k$. A
result of Hardy and Ramanujan \cite{hardy-rama} asserts that
$$
\sum_{k\leq n}\Omega(k)=n\log\log n+M'n+O\Big(\frac{n}{\log n}\Big),
$$
where
$$
M'=\gamma+\sum_{p}\left(\log\Big(1-p^{-1}\Big)+(p-1)^{-1}\right)\approx
1.0346538818.
$$
More related results can be found in Chapter V of \cite{s-m-c}.
The aim of this paper is to find an explicit version of this
result. We proceed by
$$
n!=\prod_{p\leq n}p^{v_p(n!)},
$$
standard factorization of $n!$ into primes. It is known that
$$
v_p(n!)=\sum_{k=1}^m\Big\lfloor\frac{n}{p^k}\Big\rfloor,
$$
where $\lfloor x\rfloor$ is the largest integer less than or equal
to $x$ (see for example \cite{nathan}) and
$m=m_{n,p}=\lfloor\frac{\log n}{\log p}\rfloor$. First, we
introduce some explicit (and neat) approximations for the
summation
$$
\Upsilon(n)=\sum_{p\leq n}v_p(n!).
$$
Then, considering
$$
\sum_{k\leq n}\Omega(k)=\Omega(n!)=\Upsilon(n),
$$
we obtain the main result as follows.\\
\textbf{Main Theorem.}
\textit{For every $n\geq 3$ we have}
$$
\left|\sum_{k\leq n}\Omega(k)-(n-1)\log\log(n-1)\right|<23(n-1).
$$
Note that one can modify above result to the following one:
$$
\left|\sum_{k\leq n}\Omega(k)-n\log\log n\right|<23n,
$$
which is an explicit version of the result of Hardy and Ramanujan.

\section{Proof of the Main Theorem}
Consider the inequality
\begin{equation}\label{vp-inqs}
\frac{n-p}{p-1}-\frac{\log n}{\log p}<v_p(n!)\leq \frac{n-1}{p-1},
\end{equation}
(see \cite{hassani-vp} for a proof). To get to the main theorem,
we need to approximate summations of the form $\sum_{p\leq n}f(p)$
with $f(p)=\frac{1}{\log p}$ and $f(p)=\frac{1}{p-1}$ (and more
generally, for a given function $f\in C^1(\mathbb{R}^+)$). To do
this, we use the reduction of a Riemann-Stieljes integral to a
finite sum \cite{apos-anal}, which allows us to get some ways to
evaluate the summation $\sum_{p\leq n}f(p)$; two of them are:\\

\begin{itemize}
  \item Using $\vartheta(x)=\sum\limits_{p\leq x}\log p$, which
ends to the approximation
$$
\sum_{p\leq n}f(p)=\int_{2^-}^n \frac{f(x)}{\log
x}d\vartheta(x)=\frac{f(n)\vartheta(n)}{\log
n}+\int_2^n\vartheta(x)\frac{d}{dx}\left(\frac{-f(x)}{\log
x}\right)dx,
$$
and it is known that for $x>1$, we have
$200\log^2x|\vartheta(x)-x|<793x$, and
$\log^4x|\vartheta(x)-x|<1717433x$ (see \cite{dusart} for more
details).
  \item Using $\pi(x)=\#\mathbb{P}\cap [2,x]$, which ends to
the approximation
$$
\sum_{p\leq
n}f(p)=f(x)\pi(x)+\int_2^n\pi(x)\frac{d}{dx}\big(-f(x)\big)dx,
$$
and we have some explicit bounds for $\pi(x)$ (again see
\cite{dusart} for lots of them). In this paper we will use the
following neat one:
\begin{equation}\label{pi-inq}
\pi(x)\leq\frac{x}{\log x}\Big(1+\frac{1.2762}{\log
x}\Big)\hspace{10mm}(x>1).
\end{equation}
\end{itemize}
Both of these methods are applicable for the summation
$\sum_{p\leq n}\frac{1}{p-1}$, while first method on the summation
$\sum_{p\leq n}\frac{1}{\log p}$ ends to some integrals hard to
approximate. Here, based on some known approximations for both of
these summations, which are obtained using the second method, we
give some neat bounds for them.
\begin{prop}\label{sum-1/p-1} For every $n\geq 3$, we have
$$
\log\log(n-1)-14<\sum_{p\leq n}\frac{1}{p-1}<\log\log(n-1)+23.
$$
\end{prop}
\begin{proof} It is known \cite{hassani-n!} that the inequality
$$
\log\log n+a+\frac{n}{(n-1)\log n}-\frac{1717433n}{(n-1)\log^5
n}<\sum_{p\leq n}\frac{1}{p-1},
$$
holds for $n\geq 2$ with $a\approx -11.86870152$. But, for $n\geq 3564183$ we have
$$
\log\log(n-1)-14<\log\log n+a+\frac{n}{(n-1)\log
n}-\frac{1717433n}{(n-1)\log^5 n}.
$$
Thus, for $n\geq 3564183$ we obtain
$$
\log\log(n-1)-14<\sum_{p\leq n}\frac{1}{p-1},
$$
which is also true for $2\leq n\leq 3564182$, since for these
values of $n$ the left hand side of the inequality is positive
while the right hand side is negative. Also, we have \cite{hassani-n!} the inequality
$$
\sum_{p\leq n}\frac{1}{p-1}<\log\log(n-1)+b+\frac{n}{(n-1)\log
n}+\frac{1717433 n}{(n-1)\log^5 n},
$$
for $n\geq 2$ with $b\approx 21.18095291$. On the other hand, for $n\geq
7126157$ we have
$$
b+\frac{n}{(n-1)\log n}+\frac{1717433 n}{(n-1)\log^5 n}<23.
$$
So, for $n\geq 7126157$ we obtain
$$
\sum_{p\leq n}\frac{1}{p-1}<\log\log(n-1)+23.
$$
To verify this inequality for $3\leq n\leq 7126156$, we note that
because for $p_1\leq n<p_2$ where $p_1$ and $p_2$ are two successive
primes, the left hand side is constant, while the right hand side is
increasing, therefore we only need to check this inequality for $n$
equals to prime numbers. Appendix includes the Matlab program of
doing this. The proof is completed.
\end{proof}

\begin{prop}\label{sum-1/logp} For every $n\geq 2$, we have
$$
\left|\sum_{p\leq n}\frac{1}{\log p}-\left\{\frac{n}{\log^2
n}+\frac{2n}{\log^3 n}+\frac{6n}{\log^4
n}\right\}\right|<271382\frac{n}{\log^5 n}.
$$
\end{prop}
\begin{proof} In a similar process \cite{hassani-n!}, we have
$$
\frac{n}{\log^2 n}+\frac{2n}{\log^3 n}+\frac{6n}{\log^4
n}+\frac{1607n}{100\log^5 n}-\frac{1717433n}{\log^6 n}+a<\sum_{p\leq
n}\frac{1}{\log p}\hspace{10mm}(n\geq 564),
$$
where $a\approx -16.42613005$. Also, we have
$$
\sum_{p\leq n}\frac{1}{\log p}<\frac{n}{\log^2 n}+\frac{2n}{\log^3
n}+\frac{6n}{\log^4 n}+\frac{54281n}{800\log^5
n}+\frac{1717433n}{\log^6 n}+b\hspace{10mm}(n\geq 2),
$$
where $b\approx 30.52238614$. Computation gives
$$
\frac{-271382n}{\log^5 n}<\frac{1607n}{100\log^5
n}-\frac{1717433n}{\log^6 n}+a\hspace{10mm}(n\geq 564).
$$
Also
$$
\frac{54281n}{800\log^5 n}+\frac{1717433n}{\log^6
n}+b<\frac{271382n}{\log^5 n}\hspace{10mm}(n\geq 569).
$$
Therefore, we obtain the following inequality:
$$
\left|\sum_{p\leq n}\frac{1}{\log p}-\left\{\frac{n}{\log^2
n}+\frac{2n}{\log^3 n}+\frac{6n}{\log^4
n}\right\}\right|<271382\frac{n}{\log^5 n}\hspace{10mm}(n\geq 569).
$$
A computer program verifies the above inequality for $2\leq n\leq
568$, too. The proof is complete.
\end{proof}
\hspace{-4.5mm}\textit{Proof of the Main Theorem.} Considering the
right hand side of (\ref{vp-inqs}) and the Proposition
\ref{sum-1/p-1}, for every $n\geq 3$ we have
$$
\Upsilon(n)\leq (n-1)\sum_{p\leq
n}\frac{1}{p-1}<(n-1)\log\log(n-1)+23(n-1).
$$
On the other hand, considering the left hand side of
(\ref{vp-inqs}) and the Proposition \ref{sum-1/p-1}, for every
$n\geq 3$ we have
$$
(n-1)\log\log(n-1)-14(n-1)-\mathcal{R}(n)<(n-1)\sum_{p\leq
n}\frac{1}{p-1}-\pi(n)-\log n\sum_{p\leq n}\frac{1}{\log
p}<\Upsilon(n),
$$
where
$$
\mathcal{R}(n)=\pi(n)+\log n\sum_{p\leq n}\frac{1}{\log p},
$$
and considering (\ref{pi-inq}) and the Proposition \ref{sum-1/logp},
we have
$$
\mathcal{R}(n)\leq\frac{n}{\log n}\Big(1+\frac{1.2762}{\log
n}\Big)+\frac{n}{\log n}+\frac{2n}{\log^2 n}+\frac{6n}{\log^3
n}+\frac{271382n}{\log^4 n}=\frac{2n}{\log n}+\frac{3.2762n}{\log^2
n}+\frac{6n}{\log^3 n}+\frac{271382n}{\log^4 n}.
$$
But, for $n\geq 563206$ the right hand side of this relation is
strictly less than $9(n-1)$. So, we obtain
$$
(n-1)\log\log(n-1)-23(n-1)<\Upsilon(n),
$$
for $n\geq 563206$, which holds true for $3\leq n\leq 563205$ too,
because for these values of $n$, the left hand side is positive
while the right hand side is negative. This completes the
proof.\hfill{$\Box$}

\section{Remarks for Further Studies}

\subsection{Improving the Main Result}
Of course the factor 23 in the main theorem is not optimal, and
one can improve it. But, it is the best one with our methods and
computational tools.

\subsection{Explicit Approximation of the Function $\Omega(n)$}
Concerning the main theorem, considering $n!=\Gamma(n+1)$, one can
reform the above result as
$$
\left|\Omega(\Gamma(n))-(n-2)\log\log(n-2)\right|<23(n-2),
$$
then replacing $n$ by $\Gamma^{-1}(n)$ (inverse of Gamma
function), it yields to
$$
\left|\Omega(n)-(\Gamma^{-1}(n)-2)\log\log(\Gamma^{-1}(n)-2)\right|<23\left(\Gamma^{-1}(n)-2\right).
$$
This suggests an explicit approximation for the function
$\Omega(n)$ for some special values of $n$ in terms of the inverse of Gamma function, then by
approximating $\Gamma^{-1}$, one can make it in terms of
elementary functions.

\subsection{An Extension of the Function $v_p(n!)$}
The function $v_p(n!)$, defined by
$$
n!=\prod_{p\leq n}p^{v_p(n!)},
$$
can be generalized for every positive integer $m\leq n$ instead of
prime $p\leq n$. Fix $n$ and consider canonical decomposition
$$
m=\prod_{p\leq n}p^{v_p(m)}.
$$
Same to $v_p(n!)$, we define $v_m(n!)$ in which $m^{v_m(n!)}\|n!$.
So,
$$
m^{v_m(n!)}=\prod_{p\leq n}p^{v_p(m)v_m(n!)}\Big\|\prod_{p\leq
n}p^{v_p(n!)}.
$$
Therefore, we must have $v_p(m)v_m(n!)\leq v_p(n!)$ for every prime
$p\leq n$; that is
$$
v_m(n!)\leq\min_{\substack{
p\leq n\\
v_p(m)\neq 0}}\left\{\frac{v_p(n!)}{v_p(m)}\right\}.
$$
This leads to the following definition:\\
\textbf{Definition.} For positive integers $m,n$ with $m\leq n$, we
set
$$
v_m(n!)=\left\lfloor\min_{\substack{
p\leq n\\
v_p(m)\neq 0}}\left\{\frac{v_p(n!)}{v_p(m)}\right\}\right\rfloor.
$$
Note that in the above definition, $v_p(N)$ for a positive integer
$N$ and prime $p$, is a well defined notation for the greatest
power of
$p$ dividing $N$. Related by this generalization, the following question arise to mind:\\
\textbf{Question.} Find the function $\mathfrak{F}(n)$ such that
$$
\sum_{m=1}^n v_m(n!)=\mathfrak{F}(n)\sum_{p\leq n}v_p(n!).
$$

\subsection*{Acknowledgment} We would like to express our gratitude to the referee for valuable comments.

\bigskip
\hrule\hrule\hrule
\bigskip

\begin{center}
\textsc{Appendix. Matlab program of verifying the inequality
$\sum_{p\leq n}\frac{1}{p-1}<\log\log(n-1)+23$ for prime values of
$n$}
\end{center}
{\Small\texttt{n=8000000;\\
r=primes(n);\\
s(1)=0;\\
for i=2: length(r)\\
s(i)=s(i-1)+1/(r(i)-1);\\
end\\
plot(r,s,'.',r,log(log(r)))+23,'.')}}\\\\
Final step of program plots both sides of the inequality for
comparison.

\bigskip
\hrule\hrule\hrule
\bigskip

\end{document}